\documentclass{amsart}
\usepackage[utf8]{inputenc}
\usepackage[english]{babel}
\usepackage{csquotes}
\usepackage[bookmarksdepth=3,colorlinks=true,urlcolor=blue,linkcolor=blue,citecolor=blue]{hyperref} 
\usepackage[capitalize]{cleveref}

\newcommand{\nifty}[1]{\ensuremath{\left( \infty, #1 \right)}}

\newcommand{\core}[1]{{#1}^{\text{core}}}
\newcommand{\pith}[1]{{#1}^{\text{pith}}}
\newcommand{\ipith}[1]{{#1}^{\text{pith}}}

\newcommand{\iSpc}{\icat{S}}
\newcommand{\iGrpd}{\icat{S}}

\usepackage[normalem]{ulem}

\newcommand{\catnconn}[1]{\textbf{cat-}\mathbf{#1}\textbf{-conn}}
\newcommand{\essncat}[1]{\textbf{ess-}\mathbf{#1}\textbf{-cat}}

\newcommand{\tto}{\Rightarrow}
\newcommand{\into}{\hookrightarrow}

\newcommand{\unit}{\eta}

\newcommand{\op}{^{\mathrm{op}}}
\renewcommand{\over}[1]{_{/ {#1}}}
\newcommand{\under}[1]{_{{#1}/ }}

\newcommand{\meta}[1]{\mathsf{#1}}
\newcommand{\cat}[1]{\mathit{#1}} 
\newcommand{\ccat}[1]{\underline{\mathrm{#1}}} 
\newcommand{\icat}[1]{\mathcal{#1}} 
\newcommand{\ocat}[1]{\mathbb{#1}} 
\newcommand{\so}[1]{\mathbb{#1}_{\bullet}} 
\newcommand{\ispace}[1]{\mathbf{#1}} 
\newcommand{\Hom}{\meta{Hom}} 
\newcommand{\Map}{\ispace{Map}} 
\newcommand{\iMap}{\ispace{Map}} 

\newcommand{\vMAP}[1]{\llbracket #1 \rrbracket} 
\newcommand{\MAP}{\icat{MAP}} 

\newcommand{\Grpd}{\cat{Grpd}} 

\newcommand{\grpdEmbedding}[1]{\widehat{#1}}
\newcommand{\R}{\ensuremath{\mathcal{R} }}
\newcommand{\RS}{\ensuremath{\mathcal{R}_s }}
\newcommand{\Sp}{\ensuremath{\mathrm{Sp} }}
\newcommand{\SpS}{\ensuremath{\mathrm{Sp}_S }}
\newcommand{\Ring}{\cat{cRing}}
\newcommand{\RingopS}{(\cat{cRing}\op)\over{S}}
\newcommand{\RingopR}{(\cat{cRing}\op)\over{R}}

\newcommand{\Prof}{\cat{Prof}}
\newcommand{\ProfS}{(\cat{Prof})\over{\Sp(S)}}

\newcommand{\Kalg}{K\text{-}\cat{Alg}}
\newcommand{\fKalg}{K\text{-}\cat{Alg}_{\text{fin}}}

\newcommand{\Ralg}{R\text{-}\cat{Alg}}
\newcommand{\Salg}{S\text{-}\cat{Alg}}
\newcommand{\Ralgop}{R\text{-}\cat{Alg}\op}
\newcommand{\Salgop}{S\text{-}\cat{Alg}\op}
\newcommand{\Fun}{\icat{Fun}} 
\newcommand{\lFun}{\cat{Fun}} 
\newcommand{\iFun}{\icat{Fun}} 
\newcommand{\Yo}{\icat{Y}} 
\newcommand{\dYon}{\widehat{\Yo}} 

\newcommand{\Tw}{\icat{Tw}} 

\newcommand{\N}{\mathscr{N}} 
\newcommand{\Ho}{\cat{h}} 
\newcommand{\Nerve}{\N} 
\newcommand{\icore}[1]{\left( #1 \right)^{\simeq}} 
\newcommand{\ho}[1]{\Ho #1 } 

\newcommand{\A}{\ensuremath{\icat{A}} }
\newcommand{\B}{\ensuremath{\icat{B}} }
\newcommand{\C}{\ensuremath{\icat{C}} }

\newcommand{\D}{\ensuremath{\icat{D}} }

\newcommand{\comp}{\circ }

\newcommand{\nat}{\Rightarrow }

\newcommand{\Cat}{\cat{Cat}_1} 
\newcommand{\CCat}{\ccat{Cat}} 
\newcommand{\iCat}{\icat{Cat}_{\infty}}  
\newcommand{\nCat}{\cat{Cat}}  
\newcommand{\iSpaces}{\icat{S}} 

\newcommand{\ZZ}{\ensuremath{\mathbb{Z}} }

\newcommand{\GG}{\ocat{G} }

\newcommand{\PP}{\ocat{P} }

\newcommand{\FF}{\ocat{F} }
\newcommand{\CC}{\ocat{C} }
\newcommand{\DD}{\ocat{D} }
\newcommand{\topos}[1]{\ensuremath{\icat{#1}}}

\renewcommand{\L}{\topos{L} }
\newcommand{\infcat}{\( (\infty,1)\)-category\ } 

\newcommand{\iGroth}{\icat{Groth}_\infty }

\newcommand{\DDelta}{\mathbf{\Delta} }
\newcommand{\Simplex}{\cat{\Delta} }

\newcommand{\Simp}{\Simplex}

\newcommand{\Sets}{\cat{Set} }

\newcommand{\Set}{\Sets}

\renewcommand{\S}{\icat{S} }

\newcommand{\Horn}{\Lambda }
\newcommand{\del}{\partial}
\newcommand{\isisom}{\cong}

\newcommand{\tensor}{\otimes}
\newcommand{\isequiv}{\simeq} 

\newcommand{\lSplit}{\cat{Split}}

\newcommand{\Split}{\icat{Split}}

\newcommand{\Colim}{\varinjlim}

\newcommand{\Lim}{\varprojlim}
\renewcommand{\lim}[1]{\Lim \left( #1 \right)}

\newcommand{\porth}{\ {\begin{sideways}$\!\Vdash$\end{sideways}}\ }

\newlength\squareheight
\newcommand\restr[2]{{
  \left.\kern-\nulldelimiterspace 
  #1 
  \vphantom{\big|} 
  \right|_{#2} 
  }}

\usepackage{scalerel}
\newcommand{\pb}{\scaleobj{2}{\lrcorner}}
\newcommand{\Lrcorner}{\pb}

\usepackage{geometry} 
\usepackage{bibentry} 
\usepackage{stmaryrd}
\usepackage{rotating}
\usepackage{amsmath}
\usepackage{amsthm}
\usepackage{tabularx}
\usepackage{tikz} 
\usepackage{tikz-cd} 
\usepackage{quiver} 
\usepackage{caption}
\usepackage{subcaption}
\usepackage{wrapfig}
\usepackage{multicol}

\usepackage[utf8]{inputenc}
\usepackage[cal=esstix,
scr=rsfs
]{mathalpha}
\usepackage{mathabx}

\usepackage{amsfonts} 

\usepackage{pdfsync}

\usepackage{booktabs} 
\usepackage{lipsum}
\usepackage{enumitem}

\usepackage[shortcuts]{extdash}

\newcommand{\de}[1]{\textbf{#1}}

\newtheorem{thm}{Theorem}[section]
\newtheorem{theorem}[thm]{Theorem}

\newtheorem{lemma}[thm]{Lemma}
\newtheorem{remark}[thm]{Remark}

\newtheorem{corollary}[thm]{Corollary}
\newtheorem{prop}[thm]{Proposition}

\newtheorem{definition}[thm]{Definition}
\newtheorem{principle}[thm]{Principle}

\newtheorem{xample}[thm]{Example}

\newtheorem*{thm*}{Theorem}
\newtheorem*{theorem*}{Theorem}
\newtheorem*{lemma*}{Lemma}

\usepackage{subfiles}

\usepackage[backend=biber,style=alphabetic,citestyle=alphabetic,sorting=none,maxnames=99]{biblatex}
\usepackage{todonotes}
\newcommand{\joseph}[2][]{
  \ifmmode
   \textrm{\todo[author=Joseph,backgroundcolor=green,#1]{#2}}
  \else
  \todo[author=Joseph,backgroundcolor=green,inline,#1]{#2}
  \fi
}

\newcommand{\charles}[2][]{
  \ifmmode
    \textrm{\todo[author=Charles,backgroundcolor=orange,inline,#1]{#2}}
  \else
    \todo[author=Charles,backgroundcolor=orange,inline,#1]{#2}
  \fi
}
\newcommand{\vesna}[2][]{
  \ifmmode
    \textrm{\todo[author=Vesna,backgroundcolor=orange,inline,#1]{#2}}
  \else
    \todo[author=Vesna,backgroundcolor=orange,inline,#1]{#2}
  \fi
}
\newcommand{\jeremiah}[2][]{
  \ifmmode
    \textrm{\todo[author=Jeremiah,backgroundcolor=orange,inline,#1]{#2}}
  \else
    \todo[author=Jeremiah,backgroundcolor=orange,inline,#1]{#2}
  \fi
}
\newcommand{\matt}[2][]{
  \ifmmode
    \textrm{\todo[author=Matt,backgroundcolor=orange,inline,#1]{#2}}
  \else
    \todo[author=Matt,backgroundcolor=orange,inline,#1]{#2}
  \fi
}

\begin{document}

\begin{abstract}
  In \cite{galois-theories} they prove a `Categorical Galois Theorem' for ordinary categories, and establish the main result of \cite{joyal-tierney}, along with the classical Galois theory of Rings, as instances of this more general result. 
  The main result of the present work is to refine this to a `Quasicategorical Galois Theorem', by drawing heavily on the foundation laid in \cite{kerodon}. 
  More importantly, the argument used to prove the result is intended to highlight a deep connection between factorization systems (specifically the `lex modalities' of \cite{abfj-sheaves}),  higher-categorical Galois Theorems, and Galois theories internal to higher toposes. 
  This is the first part in a series of works, intended merely to motivate the lens and prove Theorem~\ref{thm:qcat-galois}. In future work, we will delve into a generalization of the argument, and offer tools for producing applications.
\end{abstract}
\title[A Quasicategorical Galois Theorem]{The Unreasonable Efficacy of the Lifting Condition in Higher Categorical Galois Theory I: a Quasicategorical Galois Theorem}
\author{Joseph Rennie}
\maketitle


    \hypersetup{linkcolor=black}  

\section*{List of Symbols}
Throughout this paper, we assume the same framework, foundations, and notation as \cite{kerodon} except where otherwise stated. The following is a brief summary of the notation used here, which may differ from that of \cite{kerodon}.
\begin{enumerate}
\item[$\S$] $(\infty,1)$-category of Kan complexes
  \item[$\iCat$] $(\infty,1)$-category of small $(\infty,1)$-categories 
  \item[$\Cat$] $1$-category of small $1$-categories 
\item[$\Hom_{\cat{C}}(X,Y)$] Set of arrows from $X$ to $Y$ in a $1$-category
\item[$\Map_{\icat{C}}(X,Y)$] $\infty$-groupoid of arrows from $X$ to $Y$ in an $\infty$-category
\item[$\vMAP{X,Y}_{\icat{C}}$] $(\infty,1)$-category of arrows from $X$ to $Y$ in an $(\infty,1)$-category $\icat{C}$ with enhanced mapping spaces \footnote{denoted $\MAP_{\icat{C}}(X,Y)$ in \cite{ghn}.}
  \item[$\lFun(\cat{C},\cat{D})$]  $1$-category of functors from $\cat{C}$ to $\cat{D}$ in $\Cat$
  \item[$\Fun(\icat{C},\icat{D})$]  $\infty$-category of arrows from $\C$ to $\D$ in $\iCat$
    \item[$\Grpd(\C)$] $\infty$-category of groupoid objects in $\C$
\end{enumerate}

    \setcounter{section}{-1}
\section{Introduction}
\label{sec:preface}
\label{sec:introduction}
\label{sec:intro}

\subsection{History and Motivation}
One of the earliest, and perhaps most striking, demonstrations of the advantages of Category Theory as more than a mere convenience of language\footnote{which the author believes to be sufficient justification for study} came from the Galois Theory of Grothendieck established in the well-known "SGA" I-IV.
One major result of the aforementioned can be summarized as saying that certain connected, atomic Grothendieck toposes (as all those associated to separable field extensions) are each equivalent to the classifying topos of their associated Galois group. 
In \cite{joyal-tierney}, the authors extend the result of Grothendieck to say that every Grothendieck topos is equivalent to the classifying topos of a localic groupoid.
One major result of \cite{galois-theories} (the `categorical Galois theorem') then established not only a generalization of the above, but a categorification of the concepts one often associates with a proper Galois theory: splitting extensions and an appropriate notion of Galois group from which to take subgroups. The latter result will be summarized here in Section~\ref{sec:ordinary-galois}, as it is a central motivation for the project established in this paper, and foreshadowing work in its sequels.

Meanwhile, higher category theory was being developed as a tool to allow for homotopy theory to be done through a category-theoretic lens, which culminated in the definition of a higher topos in \cite{model-toposes}, with an extensive foundation laid out in \cite{htt}.
There followed many instances of higher categorical takes on various Galois-type results, including \cite{rognes-galois}, \cite{banerjee-galois}, \cite{matthew-galois}, and \cite{hoyois-galois}.
While the author of the present work was a Ph.D. student (working under the author of \cite{model-toposes}) the original project was to extend the classification of Grothendieck toposes, proven in \cite{joyal-tierney}, to higher toposes. 
Along the way, this goal gave way to that of generalizing \cite{galois-theories} to a higher categorical Galois theorem.
This paper proves such a generalization (Theorem~\ref{thm:qcat-galois-theorem}) but the method of doing so is more important. For one, it highlights (via recent results from \cite{abfj-sheaves} and its sequels) a connection to the original goal of establishing a Galois theory for higher toposes. Secondly, the framework of the proof bolsters the thesis of the current work, and the motivation for the project established, both of which we summarize in the next subsection.

\subsection{Formal Backdrop and Thesis}
We begin with the following model-independent summary of the higher categorical landscape, discussed thoroughly in \cite{kerodon}, for example.
\begin{figure}
\[\begin{tikzcd}[ampersand replacement=\&]
	\Set \\
	\Grpd \&\& {\nCat_{(1,1)}} \\
	\\
	{\Grpd_2} \&\& {\nCat_{(2,1)}} \&\& {\nCat_{(2,2)}} \\
	\& \vdots \&\& \vdots \\
	\iGrpd \&\& \iCat \&\& {\nCat_{(\infty,2)}} \& \cdots
	\arrow[""{name=0, anchor=center, inner sep=0}, "{\iota_0}"', curve={height=6pt}, hook, from=1-1, to=2-1]
	\arrow[""{name=1, anchor=center, inner sep=0}, "{\pi_0}"', curve={height=6pt}, from=2-1, to=1-1]
	\arrow[""{name=2, anchor=center, inner sep=0}, "i"', curve={height=6pt}, hook, from=2-1, to=2-3]
	\arrow[""{name=3, anchor=center, inner sep=0}, "\iota"', curve={height=12pt}, hook, from=2-1, to=4-1]
	\arrow[""{name=4, anchor=center, inner sep=0}, "{\icore{\bullet}}"', curve={height=6pt}, from=2-3, to=2-1]
	\arrow[""{name=5, anchor=center, inner sep=0}, "{\Nerve_{\bullet}}"', curve={height=12pt}, hook, from=2-3, to=4-3]
	\arrow[""{name=6, anchor=center, inner sep=0}, "{\tau_\bullet}"', curve={height=12pt}, from=4-1, to=2-1]
	\arrow[""{name=7, anchor=center, inner sep=0}, "i"', curve={height=6pt}, hook, from=4-1, to=4-3]
	\arrow[""{name=8, anchor=center, inner sep=0}, "\iota"', curve={height=12pt}, hook, from=4-1, to=6-1]
	\arrow[""{name=9, anchor=center, inner sep=0}, "\Ho"', curve={height=12pt}, from=4-3, to=2-3]
	\arrow[""{name=10, anchor=center, inner sep=0}, "{\icore{\bullet}}"', curve={height=6pt}, from=4-3, to=4-1]
	\arrow[""{name=11, anchor=center, inner sep=0}, "i"', curve={height=6pt}, hook, from=4-3, to=4-5]
	\arrow[""{name=12, anchor=center, inner sep=0}, "{\Nerve_\bullet}"', curve={height=12pt}, hook, from=4-3, to=6-3]
	\arrow[""{name=13, anchor=center, inner sep=0}, "{\pith{\bullet}}"', curve={height=6pt}, from=4-5, to=4-3]
	\arrow[""{name=14, anchor=center, inner sep=0}, "{\Nerve_\bullet}"', curve={height=12pt}, hook, from=4-5, to=6-5]
	\arrow[""{name=15, anchor=center, inner sep=0}, "{\tau_\bullet}"', curve={height=12pt}, from=6-1, to=4-1]
	\arrow[""{name=16, anchor=center, inner sep=0}, "i"', curve={height=6pt}, hook, from=6-1, to=6-3]
	\arrow[""{name=17, anchor=center, inner sep=0}, "\Ho"', curve={height=12pt}, from=6-3, to=4-3]
	\arrow[""{name=18, anchor=center, inner sep=0}, "{\icore{\bullet}}"', curve={height=6pt}, from=6-3, to=6-1]
	\arrow[""{name=19, anchor=center, inner sep=0}, "i"', curve={height=6pt}, hook, from=6-3, to=6-5]
	\arrow[""{name=20, anchor=center, inner sep=0}, "{\Ho_\bullet}"', curve={height=12pt}, from=6-5, to=4-5]
	\arrow[""{name=21, anchor=center, inner sep=0}, "{\ipith{\bullet}}"', curve={height=6pt}, from=6-5, to=6-3]
	\arrow["\dashv"{anchor=center, rotate=-180}, draw=none, from=1, to=0]
	\arrow["\dashv"{anchor=center, rotate=-90}, draw=none, from=4, to=2]
	\arrow["\dashv"{anchor=center, rotate=-180}, draw=none, from=6, to=3]
	\arrow["\dashv"{anchor=center, rotate=-180}, draw=none, from=9, to=5]
	\arrow["\dashv"{anchor=center, rotate=-90}, draw=none, from=10, to=7]
	\arrow["\dashv"{anchor=center, rotate=-90}, draw=none, from=13, to=11]
	\arrow["\dashv"{anchor=center, rotate=-180}, draw=none, from=15, to=8]
	\arrow["\dashv"{anchor=center, rotate=-90}, draw=none, from=18, to=16]
	\arrow["\dashv"{anchor=center, rotate=-180}, draw=none, from=17, to=12]
	\arrow["\dashv"{anchor=center, rotate=-180}, draw=none, from=20, to=14]
	\arrow["\dashv"{anchor=center, rotate=-90}, draw=none, from=21, to=19]
\end{tikzcd}\]
\caption{}
\label{fig:ncat-context}
\end{figure}
Note that the above can be made precise\footnote{Full precision demands that we fix a model, and specify the bullets in the vertical functors above with the appropriate construction. 
In the present work, all of the relevant functors above will come from \cite{kerodon}, but it will always be clear from context which functor is to be used. For a detailed comparison of different (Quillen equivalent) choices of models for the above, we refer the interested reader to \cite{bergner-rezk-survey}, \cite{bergner-rezk-i}, and \cite{bergner-rezk-ii}.}
\footnote{Full precision also requires that we specify what kind of category each of the above is. Since we only care about the categories depicted, we can simply demand that they are all (perhaps trivially) in the very large $(\infty,2)$-category of large $(\infty,2)$-categories. 
Without specifying exactly how, we merely note that this can be done in many ways, and refer the interested reader to \cite{size-stuff} for a thorough (and insightful) exposition of this.}, but our intent is only to highlight its conceptual value.
Specifically, the chain of fully faithful right adjoints will be shown to interact nicely with the given proof (centered around lifting conditions) of the quasicategorical Galois theorem presented here. 
This is one of the main theses of the current work, and foreshadows future work in which we extend the main result proven here using only the formal backdrop of the above diagram.
Indeed, this work can be seen as one such formal extension along the vertical adjunctions of the first column above.

At this point it is crucial to emphasize the other central aims of this project:
\begin{enumerate}
  
  \item In \cite{abfj-sheaves}, they work out a similar interaction of certain factorization systems with localizations of $(\infty,1)$-toposes, and each of the above can be thought of as the ur-example of an $(\infty, n)$-topos.  Thus, one can imagine (as investigated in the successors to the current work) establishing an internal galois theorem for arbitrary $(\infty, n)$-toposes which yields an analogue to the main result here for any $(\infty,n-1)$-category. 
  \item The categorical Galois theorem of \cite{galois-theories} specifies down to the Galois theory internal to any ordinary (Grothendieck) topos established in \cite{joyal-tierney}.
\item Finally, the proof given here suggests the existence of a `model-independent' proof of a higher categorical Galois theorem. Recent work of \cite{riehl-verity} and \cite{nima-model-independent} establishes a formal meaning of the phrase above which we use in future work to investigate the desired model-independent proof.
\end{enumerate}

The purpose of this paper is to summarize the conceptual arc of \cite{galois-theories}, and prove the quasicategorical analogue of their main result in a way that motivates and establishes a conceptual underpinning of the above project to be carried out in sequels to the present work.
\section{Ordinary Categorical Galois Theories}
\label{sec:categorical-galois-theories}
\label{sec:ordinary-galois}

This section reviews the specification of the ordinary categorical Galois theorem of \cite{galois-theories} to give the Galois theory of rings as an instance.
The purpose is threefold:
\begin{enumerate}
  \item it summarizes the essentials of a justification of the `Galois' in `[ordinary] categorical Galois Theorem' to the uninitiated,
  \item it foreshadows a general method of producing such results as formal consequences of certain adjunctions (which we refine to the higher categorical setting in future work), and
  \item it motivates our definitions used in the quasicategorical Galois theorem. 
\end{enumerate}
\label{ssec:galois-theories}
Classically, the Galois theory of field extensions is stated as an anit-equivalence between two posets associated to a fixed Galois field extension: 1) the poset of field sub-extensions and 2) the poset of subgroups of the group of automorphisms which fix the base field. 
We really start to see a Galois \textit{principle} (or metatheorem) emerge when we compare this to the theory of covering spaces.
Here we again have an equivalence between two posets associated to a regular covering of spaces: 1) the poset of subcoverings, and 2) the poset of subgroups of the group of deck transformations.

Seeking to unify these two results leaves us with two concepts to categorify:
\begin{enumerate}
  \item splitting extensions, and
  \item subgroups of the group of automorphisms fixing a subobject.
\end{enumerate}
Additionally, we hope to do so in a way that makes the following principal a mathematically explicit one.
\begin{principle}
  The category of ``splitting objects'' within a given extension $\sigma:K \to L$ is \\ anti-equivalent to the category of ``subgroups of the Galois Group(oid).''
\end{principle}
The ordinary category formulation of the above principal is a large part of \cite{galois-theories}. 
Its quasicategorical analogue is the main result (but not the main purpose) of this paper.
While everything in this section comes from \cite{galois-theories}, we have adapted the language to more easily merge into the higher categorical project established here.
The definitions we use in the higher categorical program are motivated by their analogues in \cite{galois-theories}.

\subsection{Splitting Field Extensions}
\label{ssec:splitting-fields}
The notion of a splitting object for a field extension $K\to L$ goes by the name of ``Galois sub-extension'', which we recall here.
Typically, this property is given as a conjunction of two other  properties of field extensions: normality and separability.
Our aim here is to express these properties as equations, and thus facilitate the categorification of the notion of splitting object.
\begin{definition}
Given a field extension $\sigma : K\to L$ a \de{minimal polynomial} of $l\in L$ is an irreducible polynomial $p(x)\in K[x]$ satisfying $p(l)=0$. We say $sigma$ is \de{algebraic} if every element of $l$ has a minimal polynomial.
\end{definition}
\begin{definition}
A field extension $\sigma:K \to L$ is \de{normal} if it is algebraic and every minimal polynomial factors into linear polynomials in $L[x]$.
\end{definition}
\begin{definition}
A field extension $\sigma:K \to L$ is \de{separable} if it is algebraic and every minimal polynomial factors into distinct (linear) factors in $\overline{K}[x]$ where $\overline{K}$ is the algebraic closure of $K$.
\end{definition}
Here we might simply define a Galois extension as a normal and separable one, leading us to internalize the notion using two equations.
However, the following results will allow us to simplify this definition to one equation.
  \begin{lemma}[{\cite[Ch 2]{galois-theories}}]
    Let $\fKalg$ denote the category of finite-dimensional $K$-algebras. Then each of the following holds:
      \[ L\tensor_K \frac{K\left[ x \right]}{\langle p(x) \rangle}
          \isisom \frac{L\left[ x \right]}{\langle p(x) \rangle}\]
     \[ L\text{-roots of } p(x)\in K[x] \isisom Hom_{\fKalg} (\frac{K[x]}{p(x)}, L)\]
  \end{lemma}
  Here we note that the definition of a Galois extension is entirely centered around polynomials, their roots, and irreducibility.
  By the above, each of these concepts makes sense in the more general context of $K$-algebras other than $L$.
In the case where the $K$-algebra is a field extension, separablility is expressed in the isomorphism labelled (1) below, and normality is expressed in the isomorphism labelled (2).
\begin{lemma}[{\cite[Ch 2]{galois-theories} }]
  \label{lem:galois-def-fields}
If $\sigma: K \to L$ is a Galois extension and $a\in L$ with minimal polynomial $p_a(x)\in K[x]$,  then we have
    \[L^{Hom_{\Kalg}(\frac{K[x]}{p_a(x)},L)} \isisom L^{\text{\# of roots}} \isisom_{(1)} L^{deg(p)} \isisom_{(2)} \frac{L[x]}{p_a(x)} \isisom L\tensor_K \frac{K[x]}{p_a(x)}\]
\end{lemma}
In the next section we rely on the recognition of the above as a natural isomorphism coming from an adjunction to arrive at a more general notion of splitting extensions applicable to other catgories in place of the category of finite fields.
\subsection{Splitting Objects}
\label{ssec:splitting-objects}
\textit{(This section is essentially a rephrasing of a few chapters of \cite{galois-theories} to be more suited to adaptation to the higher categorical context. Along the way, we show how certain adjunctions naturally lead to applications of both Theorem~\ref{thm:categorical-galois-theorem} and Theorem~\ref{thm:quasicategorical-galois-theorem}.)}

There is a well-known principle, deep and pervasive in much of professional Mathematics through the last century, which says that algebra is ``dual'' to geometry.
Roughly speaking,  we often find that results obtained studying a collection of algebraic objects automatically yield results about some corresponding collection of geometric ones.
Often this can be made precise as an adjunction between a category of geometric objects and the opposite of its algebraic counterpart.
Here, we begin with one such example and show how it naturally leads us to a Galois theory for (commutative unital) Rings.
In doing so, we arrive at an internalization of ``splitting object'' which immediately reveals the isomorphism of Lemma~\ref{lem:galois-def-fields} as an instance.

The duality of interest to us in this exposition comes from the Stone Duality (between boolean algebras and profinite spaces) being restricted along the inclusion of commutative rings into boolean algebras.
\begin{theorem}[4.3.2 of \cite{galois-theories}]
Let $\Ring$ denote the category of commutative unital rings, and $\Prof$ denote the category of profinite topological spaces.
    We have an adjunction:
\[\begin{tikzcd}[ampersand replacement=\&]
	\Ring\op \& \Prof
	\arrow[""{name=0, anchor=center, inner sep=0}, "\Sp", shift left=2, from=1-1, to=1-2]
	\arrow[""{name=1, anchor=center, inner sep=0}, "\R", shift left=2, from=1-2, to=1-1]
	\arrow["\dashv"{anchor=center, rotate=-90}, draw=none, from=0, to=1]
\end{tikzcd}\]
where $\Sp$ is the ``Pierce Spectrum'' and $\R$ is the functor defined on objects by sending a profinite space $X$ to the ring of continuous functions $\C(X,\ZZ)$ (where $\ZZ$ has the discrete topology).
  \end{theorem}
  Note that the precise definitions of $\Prof$, $\Sp$, and $\R$ above are irrelevant for our purposes.
  What we will see is that formal properties of the adjunction will suffice to guide our discussion, thereby suggesting a more general framework for producing such results. Toward this aim we begin with the following observation.
  \begin{lemma}[4.3.4 of \cite{galois-theories}]
    \label{obs:ring-pushout}
The diagram below is a pushout in $\Ring$ (a pullback in $\Ring\op$)
\[\begin{tikzcd}[ampersand replacement=\&]
	A \& C \\
	B \& {B\otimes_A C}
	\arrow[from=1-1, to=2-1]
	\arrow[from=2-1, to=2-2]
	\arrow[from=1-2, to=2-2]
	\arrow[from=1-1, to=1-2]
	\arrow["\lrcorner"{anchor=center, pos=0.125, rotate=180}, draw=none, from=2-2, to=1-1]
\end{tikzcd}\]
\end{lemma}
More generally, whenever we have two categories $\C$ and $\D$ along with an adjunction
\[\begin{tikzcd}[ampersand replacement=\&]
	\C \& \D
	\arrow[""{name=0, anchor=center, inner sep=0}, "\L", shift left=2, from=1-1, to=1-2]
	\arrow[""{name=1, anchor=center, inner sep=0}, "\R", shift left=2, from=1-2, to=1-1]
	\arrow["\dashv"{anchor=center, rotate=-90}, draw=none, from=0, to=1]
\end{tikzcd}\]
where $\C$ has pullbacks, we automatically get an adjunction (as below) for each object $C\in \C$.
\[\begin{tikzcd}[ampersand replacement=\&]
	\C\over{C} \& \D\over{\L(C)}
	\arrow[""{name=0, anchor=center, inner sep=0}, "\L_C", shift left=2, from=1-1, to=1-2]
	\arrow[""{name=1, anchor=center, inner sep=0}, "\R_C", shift left=2, from=1-2, to=1-1]
	\arrow["\dashv"{anchor=center, rotate=-90}, draw=none, from=0, to=1]
\end{tikzcd}\]
Thus, Lemma~\ref{obs:ring-pushout} says we have an adjunction (as below) for any $S\in \Ring$.
\[\begin{tikzcd}[ampersand replacement=\&]
	\RingopS \& \ProfS
	\arrow[""{name=0, anchor=center, inner sep=0}, "\SpS", shift left=2, from=1-1, to=1-2]
	\arrow[""{name=1, anchor=center, inner sep=0}, "\RS", shift left=2, from=1-2, to=1-1]
	\arrow["\dashv"{anchor=center, rotate=-90}, draw=none, from=0, to=1]
\end{tikzcd}\]
We note that for $R\in \Ring$ we have $\Ralg \isequiv \Ring\under{R}$, and thus the above adjunction can be rewritten to more closely resemble the context of the Galois theory of fields at the end of the previous section.
\[\begin{tikzcd}[ampersand replacement=\&]
	\Salgop \& \ProfS
	\arrow[""{name=0, anchor=center, inner sep=0}, "\SpS", shift left=2, from=1-1, to=1-2]
	\arrow[""{name=1, anchor=center, inner sep=0}, "\RS", shift left=2, from=1-2, to=1-1]
	\arrow["\dashv"{anchor=center, rotate=-90}, draw=none, from=0, to=1]
\end{tikzcd}\]
The following result directly connects the isomorphism defining Galois extensions to a property easily stated in terms of the above adjunction.
\begin{lemma}[{\cite{gtheories}}]
  \label{fact:splitting-unit}
  If $\sigma: K \to L$ is a Galois finite field extension, and $A\in \Kalg$ is split by $\sigma$, then
  \[  L^{Hom_K(A,L)}\isisom \R (\Sp (L\tensor_K A))\]
\end{lemma}
Combining Lemma~\ref{fact:splitting-unit} with Lemma~\ref{lem:galois-def-fields} tells us that $A$ is a Galois sub-extension of $\sigma$ if only if
  \[  L\tensor_K A\isisom \R (\Sp (L\tensor_K A)).\]
This is more succinctly expressed as saying the $(L\tensor_K A)$-component of the unit of the adjunction is an isomorphism.\footnote{This is even more succinctly, and suggestively, expressed by saying the splitting algebras of $\sigma: K \to L$ are precisely the fields whose image under $L\tensor_K \bullet$ are left fixed points of the adjunction associated to $\sigma$.}
At this point, we can define splitting algebras with respect to a given map of rings:
\begin{definition}
     A given $B\in \Salg$ \de{splits} if $\unit_B \op : \RS( \SpS (B)) \to B$ is an isomorphism in $\Salg$.

     Given $\sigma:R\to S$ in $\Ring$, we say $A\in \Ralg$ \de{is split by $\sigma$} if $S\tensor_R A$ splits in the above sense.
       We let $\lSplit_R(\sigma)$ denote the full subcategory of $\Ralg$ consisting of objects split by $\sigma$.
\end{definition}

The restricted Stone Duality which began our discussion happens to have the property that $\RS$ is a fully faithful functor, and thus the counit is a natural isomorphism (i.e. every component is an isomorphism).
We can then restrict the adjunction to the subcategory of $\Salgop$ whose objects give components of the unit which are isomorphisms.
In doing so we force the adjunction to become an equivalence of categories. Alternatively, $\RS$ being fully faithful means thta $\ProfS$ is isomorphic to a full subcategory of $\Salgop$ which defines $\lSplit_S(S)$.
Note that the diagram above is a definition of $\lSplit_S(S)$ in that the top square is a pseudolimit.
Furthermore, given $\sigma: S\to R$ in $\Ring\op$ we get a functor \[S\tensor_R \bullet : \RingopR \to \RingopS ,\]
or equivalently \[S\tensor_R \bullet : \Ralgop \to \Salgop .\]
This lets us define a notion of \de{splitting algebras relative to $\sigma$} via the square below being a pseudolimit.
\[\begin{tikzcd}[ampersand replacement=\&]
	\& {\lSplit_R(\sigma)} \& \Ralgop \\
	\ProfS \& {\lSplit_S(S)\op} \& \Salgop
	\arrow[hook, from=1-2, to=1-3]
	\arrow[from=1-2, to=2-2]
	\arrow["\lrcorner"{anchor=center, pos=0.125}, draw=none, from=1-2, to=2-3]
	\arrow["{S\tensor_R \bullet}", from=1-3, to=2-3]
	\arrow["\sim", tail reversed, from=2-2, to=2-1]
	\arrow[shift right=2, hook, from=2-2, to=2-3]
\end{tikzcd}\]

We finally arrive at the categorical notion of splitting objects, by following the path suggested by core results in category theory.
Namely, any category $\C$ with finite limits yields a pseudofunctor $\C\over{\bullet}:\C\op\to \CCat$, where $\CCat$ is the $2$-category of ``categories.'' 
In fact, the adjunction
\[\begin{tikzcd}[ampersand replacement=\&]
	\RingopS \& \ProfS
	\arrow[""{name=0, anchor=center, inner sep=0}, "\SpS", shift left=2, from=1-1, to=1-2]
	\arrow[""{name=1, anchor=center, inner sep=0}, "\RS", shift left=2, from=1-2, to=1-1]
	\arrow["\dashv"{anchor=center, rotate=-90}, draw=none, from=0, to=1]
\end{tikzcd}\]
yields two such pseudofunctors into $\CCat$ along with a (pseudo)natural transformation between them:
\begin{equation}
  \label{con:splitting-objects}
  \tag{S1}
  \begin{tikzcd}[ampersand replacement=\&]
	\Ring \& S \& \Ring \& S \\
	\CCat \& \ProfS \& \CCat \& \RingopS \\
	{F_1} \& {F_2} \& S \& R
	\arrow["{F_1}", from=1-1, to=2-1]
	\arrow["{F_2}", from=1-3, to=2-3]
	\arrow[maps to, from=1-2, to=2-2]
	\arrow[maps to, from=1-4, to=2-4]
	\arrow["\R_{\bullet}", Rightarrow, from=3-1, to=3-2]
	\arrow["\sigma"', from=3-4, to=3-3]
\end{tikzcd}
\end{equation}
\begin{remark}
	In order for $\R_\bullet$ above to satisfy the necessary coherences to be a natural transformation, we rely on the fact that the original adjunction satisfied a property referred to in hte literature by a few different names. For example, \cite{townsend-proof} refers to this as `(stable) frobenius reciprocity'.
\end{remark}
We will show that with the above data, we could just as easily define $\lSplit_R(\sigma)$ by the following pseudolimit,
\[\begin{tikzcd}[ampersand replacement=\&]
	{\lSplit_R (\sigma)} \& {F_2(R)} \\
	{F_1(S)} \& {F_2(S)}
	\arrow["{\alpha_S}"', from=2-1, to=2-2]
	\arrow["{F_2(\sigma)}", from=1-2, to=2-2]
	\arrow[dashed, from=1-1, to=2-1]
	\arrow[dashed, from=1-1, to=1-2]
	\arrow["\lrcorner"{anchor=center, pos=0.125}, draw=none, from=1-1, to=2-2]
\end{tikzcd}\]
thereby motivating the following definition.
  \begin{definition}
    \label{def:splitting-object}
    Given a category $\C$, two functors $F_1,F_2:\C\op \to \CCat$, a natural transformation ${\alpha:F_1 \tto F_2}$, and an arrow $\sigma: S\to R$ in $\C$, we define $\lSplit_\alpha(\sigma)$ via the following pseudolimit in $\CCat$:
\[\begin{tikzcd}[ampersand replacement=\&]
	{\lSplit_\alpha (\sigma)} \& {F_2(R)} \\
	{F_1(S)} \& {F_2(S)}
	\arrow["{\alpha_S}"', from=2-1, to=2-2]
	\arrow["{F_2(\sigma)}", from=1-2, to=2-2]
	\arrow[dashed, from=1-1, to=2-1]
	\arrow[dashed, from=1-1, to=1-2]
	\arrow["\lrcorner"{anchor=center, pos=0.125}, draw=none, from=1-1, to=2-2]
\end{tikzcd}\]
  \end{definition}
  The following diagram shows how Definition~\ref{def:splitting-object} generalizes that of $\lSplit_R(\sigma)$.
\[\begin{tikzcd}[ampersand replacement=\&]
	{\lSplit_\alpha (\sigma)} \& {\lSplit_R(\sigma)\op} \& \Ralgop \\
	\ProfS \& {\lSplit_S(1_S)} \& \Salgop
	\arrow[dashed, from=1-1, to=1-2]
	\arrow[""{name=0, anchor=center, inner sep=0}, curve={height=-18pt}, from=1-1, to=1-3]
	\arrow[from=1-1, to=2-1]
	\arrow[hook, from=1-2, to=1-3]
	\arrow[from=1-2, to=2-2]
	\arrow["\lrcorner"{anchor=center, pos=0.125}, draw=none, from=1-2, to=2-3]
	\arrow[""{name=1, anchor=center, inner sep=0}, "{S\tensor_R \bullet}", from=1-3, to=2-3]
	\arrow["\sim"', tail reversed, from=2-2, to=2-1]
	\arrow[hook, from=2-2, to=2-3]
	\arrow["\lrcorner"{anchor=center, pos=0.125}, shift right=5, draw=none, from=0, to=1]
\end{tikzcd}\]
More explicitly, the right square is a pseudolimit by definition, so the dotted arrow exists by a defining property of the pseudolimit on the right.
The square on the left is then a pseudolimit because the outer and right squares are pseudolimits.
Since $\SpS$ restricts to an equivalence between $\lSplit_S(1_S)$ and $\ProfS$, the dotted arrow, a pullback of an equivalence, is also an equivalence.

Thus, we have a categorical notion of splitting object which subsumes the classical notion of a Galois sub-extension of rings (and thus subsumes the analogue for fields as well).
\subsection{``$\GG$-subgroups'' more generally}
Having established and motivated a categorified definition of splitting objects (namely those in $\lSplit_\alpha(\sigma)$), we now move on to do the same for the second required piece of `a reasonable Galois theorem': the collection of subgroups of a group. 
However, motivated by the classical theory of all covering spaces, not just those of connected spaces, we must widen our scope beyond groups to groupoids. 
Categorifying groupoids is much more well-trodden territory, which we summarize briefly here.

\begin{figure}[h]
     \centering
     \begin{subfigure}[b]{0.3\textwidth}
         \centering
  \[
    \begin{tikzcd}[ampersand replacement =\& ,row sep=small ,column sep = tiny]
     C_2 \ar[dr, dotted] \ar[ddr, swap,"\del_1"] \ar[drr,"\del_2"] \& \& \\
     \& \Horn_2^0(\CC) \ar[r] \ar[d] \arrow[dr, phantom, "\lrcorner", very near start]\& C_1 \ar[d,"s"] \\
     \& C_1 \ar[r, swap,"s"]\& C_0
      \end{tikzcd}
  \]
  \caption{The pullback which yields post-inverses in $\CC$.}
	     \label{fig:pullback-horns-a}
     \end{subfigure}
     \hfill
     \begin{subfigure}[b]{0.3\textwidth}
         \centering
  \[
    \begin{tikzcd}[ampersand replacement =\& ,row sep=small ,column sep = tiny]
     C_2 \ar[dr, dotted] \ar[ddr, swap,"\del_0"] \ar[drr,"\del_2"] \& \& \\
     \& \Horn_2^1(\CC) \ar[r] \ar[d] \arrow[dr, phantom, "\lrcorner", very near
     start]\& C_1 \ar[d,"t"] \\
     \& C_1 \ar[r, swap,"s"]\& C_0
      \end{tikzcd}
  \]
  \caption{The pullback which yields compositions in $\CC$.}
     \end{subfigure}
     \hfill
     \begin{subfigure}[b]{0.3\textwidth}
         \centering
  \[
    \begin{tikzcd}[ampersand replacement =\& ,row sep=small ,column sep = tiny]
     C_2 \ar[dr, dotted] \ar[ddr, swap,"\del_0"] \ar[drr,"\del_1"] \& \& \\
     \& \Horn_2^2(\CC) \ar[r] \ar[d] \arrow[dr, phantom, "\lrcorner", very near
     start]\& C_1 \ar[d,"t"] \\
     \& C_1 \ar[r, swap, "t"]\& C_0
      \end{tikzcd}
  \]
  \caption{The pullback which yields pre-inverses in $\CC$.}
     \end{subfigure}
  \\
     \begin{subfigure}[b]{0.3\textwidth}
         \centering
\[\begin{tikzcd}
	& y \\
	x && z
	\arrow["g"{description}, dotted, from=1-2, to=2-3]
	\arrow["f"{description}, from=2-1, to=1-2]
	\arrow[""{name=0, anchor=center, inner sep=0}, "h"{description}, from=2-1, to=2-3]
	\arrow["{\circlearrowright_\alpha}"{description}, draw=none, from=0, to=1-2]
\end{tikzcd}\]
\caption{a visual representation of pullback (A) above.}
       \label{fig:visual-horns-a}
     \end{subfigure}
     \hfill
     \begin{subfigure}[b]{0.3\textwidth}
         \centering
\[
	\begin{tikzcd}
	& y \\
	x && z
	\arrow["g"{description}, from=1-2, to=2-3]
	\arrow["f"{description}, from=2-1, to=1-2]
	\arrow[""{name=0, anchor=center, inner sep=0}, "h"{description}, dotted, from=2-1, to=2-3]
	\arrow["{\circlearrowright_\alpha}"{description}, draw=none, from=0, to=1-2]
\end{tikzcd}
\]
\caption{a visual representation of pullback (B) above.}
     \end{subfigure}
     \hfill
     \begin{subfigure}[b]{0.3\textwidth}
         \centering
\[\begin{tikzcd}
	& y \\
	x && z
	\arrow["g"{description}, from=1-2, to=2-3]
	\arrow["f"{description}, dotted, from=2-1, to=1-2]
	\arrow[""{name=0, anchor=center, inner sep=0}, "h"{description}, from=2-1, to=2-3]
	\arrow["{\circlearrowright_\alpha}"{description}, draw=none, from=0, to=1-2]
\end{tikzcd}\]
\caption{a visual representation of pullback (C) above.}
     \end{subfigure}
        \caption{Diagrams above internalize the Kan conditions (just on 2-cells).}
        \label{fig:visual-horn-conditions}
  \label{fig:visual-horns}
        \label{fig:horn-conditions}
  \label{fig:pullback-horns}
\end{figure}
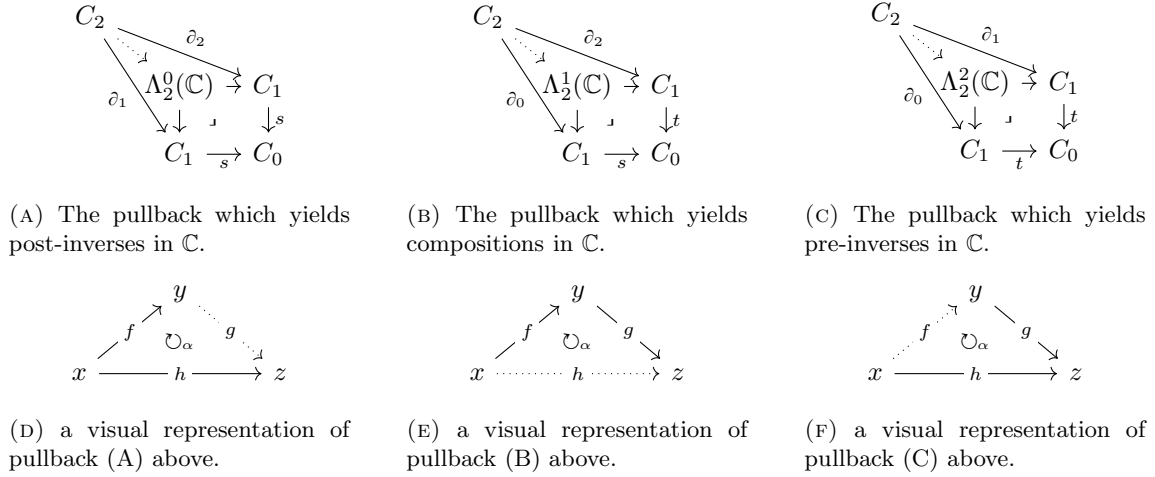

Given a discrete group $G$, every subgroup comes with an action (by translation) of $G$.
Thus we have a faithful embedding of the poset of subgroups into the category of $G$-sets.
However, even for the classical Galois theory of (possibly) infinite field extensions, it was known that the relevant groups carried the structure of a profinite topological space, and the subgroups of interest were actually only the closed subgroups. 
Here, the categorical lens simplifies the language: fixing a profinite group $G$, the poset of closed subgroups again embeds into the category of $G$-actions on profinite topological spaces.
The rest of this section is merely a recap of the meaning of all of the terms used above, with the goal of motivating the analogous definitions in the higher categorical setting.
The following arc is motivated by the conception of groupoid objects as category objects where the ``arrow-object'' consists only of invertible arrows.
\begin{definition}
Let $\C$ be a category.
  An \de{internal precategory in $\C$} is a truncated simplicial map $\CC: \Simp_{\leq 2}\to \C$ .
\end{definition}
  We can picture this as:
\[
  \begin{tikzcd}[ampersand replacement=\&]
C_2
\arrow[r,  "\del_0" ,shift left=7]
\arrow[r, leftarrow , shift left=5]
\arrow[r, "compose"]
\arrow[r, leftarrow, shift right=2]
\arrow[r, "\del_2" ,shift right=4, swap]
\& C_1
\arrow[r, "source", shift left=5]
\arrow[r, "1_{(\bullet)}", leftarrow]
\arrow[r, "target",  shift right=5]
\& C_0
\end{tikzcd}
\]
where $C_0$ is the ``object of objects'', $C_1$ is the ``object of arrows'', and $C_2$ is the ``object of compositions.''
In what follows, we use the above notation for every instance of an internal precategory, and abbreviate the source map by $s$, the target map by $t$, and `compose' by $c$.
\begin{definition}
  Given an internal precategory $\CC$ in $\C$, we say it is an \de{internal groupoid} if the
  maps required to commute by functoriality of $\CC$ are actually pullbacks in $\C$.
\end{definition}
Explicitly, the definition asks that the dotted maps of Figure~\ref{fig:horn-conditions} be isomorphisms in $\C$. The corresponding visuals as $\CC$-diagrams for each pullback in $\C$ are shown in Figure~\ref{fig:visual-horns}.
For example, $\Horn_2^0(\CC)$ is the object representing so called ``horns'' in $\CC$, so Figure~\ref{fig:pullback-horns-a} expresses the ability to fill all $\CC$-diagrams of the form given in Figure~\ref{fig:visual-horns-a} with an $\CC$-arrow $g$, and a $2$-cell $\alpha$ exhibiting the commutativity of the $\CC$-diagram.
\\  \begin{xample}
    Let $\C$ be a category with finite limits.
  \begin{itemize}
  \item Every $K\in \C$ yields a discrete internal groupoid $\GG_K$ with all maps
    being $1_K$
    \item Every morphism $\sigma:L\to K$ yields an internal groupoid $\GG_\sigma$ (called the
      kernel of $\sigma$):
      \[
        \begin{tikzcd}[ampersand replacement=\&]
L\times_K L\times_K L
\arrow[r,  shift left=4]
\arrow[r, leftarrow , shift left=2]
\arrow[r ]
\arrow[r, leftarrow, shift right=2]
\arrow[r, shift right=4]
\& L\times_K L
\arrow[r,  shift left=2]
\arrow[r,  leftarrow]
\arrow[r,   shift right=2]
\& L
        \end{tikzcd}
      \]
      \end{itemize}

  \end{xample}

In order to properly categorify the notion of a groupoid action, we rely on intuition from Algebraic Topology, where sheaves on the delooping space of a group correspond exactly with actions of the group on some space.
Thus, we seek to generalize the following definition to one which is relative to a functor $F:\C\op \to \CCat $
  \begin{definition}
    \label{def:internal-presheaf-category}
    An \de{internal presheaf} on an internal groupoid $\GG$ is a precategory $\PP$ and
    a natural transformation (internal functor) $\FF:\PP \nat \GG$
    such that all of the required commuting squares are actually pullbacks.
\[
  \begin{tikzcd}[ampersand replacement=\&]
P_2  \arrow[d, "f_2"]
\arrow[r,  shift left=4]
\arrow[r, leftarrow , shift left=2]
\arrow[r ]
\arrow[r, leftarrow, shift right=2]
\arrow[r, shift right=4]
\& P_1 \arrow[d, "f_1"]
\arrow[r,  shift left=2]
\arrow[r,  leftarrow]
\arrow[r,   shift right=2]
\& P_0 \arrow[d, "f_0"] \\
G_2
\arrow[r,  shift left=4]
\arrow[r, leftarrow , shift left=2]
\arrow[r ]
\arrow[r, leftarrow, shift right=2]
\arrow[r, shift right=4]
\& G_1
\arrow[r,  shift left=2]
\arrow[r,  leftarrow]
\arrow[r,   shift right=2]
\& G_0
\end{tikzcd}
\]

   This defines the category we denote $\C^\GG$.
  \end{definition}
  The generalization we seek is motivated by the following result.
  \begin{lemma*}[ 7.2.7 of \cite{gtheories} ]
   The category of presheaves on $\GG\in \Grpd(\C)$ is equivalent to the $2$-limit
   in $\CCat$ of the diagram below:
\[
  \begin{tikzcd}[ampersand replacement=\&]
\C/G_2
\arrow[r,  leftarrow ,"\del_0^*" ,shift left=10]
\arrow[r,  shift left=5]
\arrow[r, leftarrow ,"c^*"]
\arrow[r,  shift right=5]
\arrow[r, leftarrow ,"\del_2^*" ,shift right=10]
\& \C/G_1
\arrow[r, leftarrow, "s^*", shift left=5]
\arrow[r, "1_{(\bullet)}"]
\arrow[r, leftarrow, "t^*",  shift right=5]
\& \C/G_0
\end{tikzcd}
\]

  \end{lemma*}
  This result lets us define a more general notion of internal presheaves \textit{with respect to a functor} ${F:\C\op \to \CCat}$, which has Definition~\ref{def:internal-presheaf-category} as an instance.
  \begin{definition}
   Given a pseudo-functor $F:\C^{op}\to \CCat$ and an internal groupoid $\GG\in \Grpd(\C)$ we define the \de{category of $F$-presheaves on $\GG$}, $F^\GG$, to be the
   pseudo-limit $Lim(F\comp \GG)$.
  \end{definition}
Explicitly, using the ``slice functor'' $C\to \C\over{C}$ in place of $F$ yields Definition~\ref{def:internal-presheaf-category}.
\subsection{Categorical Galois Theorem}
\label{ssec:categorical-galois-theorem}
We recall Context~\ref{con:splitting-objects} here.
\begin{equation*}
\begin{tikzcd}[ampersand replacement=\&]
	\Ring \& S \& \Ring \& S \\
	\CCat \& \ProfS \& \CCat \& \RingopS \\
	{F_1} \& {F_2} \& S \& R
	\arrow["{F_1}", from=1-1, to=2-1]
	\arrow["{F_2}", from=1-3, to=2-3]
	\arrow[maps to, from=1-2, to=2-2]
	\arrow[maps to, from=1-4, to=2-4]
	\arrow["\alpha", Rightarrow, from=3-1, to=3-2]
	\arrow["\sigma"', from=3-4, to=3-3]
\end{tikzcd}
\end{equation*}

By placing reasonable conditions on $\sigma$  and $\alpha$, \cite{galois-theories} proves the $1$-categorical Galois theorem which states that $\lSplit_\alpha(\sigma)\isequiv F_1^{\GG_\sigma}$. 
The condition on $\alpha$ is motivated by the example context of $\Ring$. Namely, they require that $\alpha$ be fully faithful on each component.
The condition on $\sigma$ is motivated by the descent theory of modules over rings.
However, here, we do not motivate the categorification of descent, as doing so would not add to our thesis here in the way that categorifying splitting objects does.\footnote{See \cite[\href{https://stacks.math.columbia.edu/tag/023F}{Tag 023F}]{descent-modules} for a thorough exposition of this.}

\begin{lemma}[\cite{galois-theories}]
       Given a $F:\C\op \to \CCat$, and an internal functor $\FF:\CC\to \DD$ between internal groupoids yields a \de{comparison functor}, $\kappa^\FF$.
        \[
        \kappa^\FF: F^\DD \to F^\CC
        \]
\end{lemma}
        In particular, for any $\sigma:K \to L$ in $\C$, there is an internal functor \(\FF_\sigma:\GG_\sigma\to \GG_K \) given by projection (then $\sigma$) at each level.
        This internal functor yields a corresponding comparison functor, which is a central piece in the study of descent.
  \begin{definition}
    Given $F:\C^{op}\to \CCat$, we say that $\sigma:L\to K$ is of \de{$F$-effective descent} if
   for the canonical groupoid-map $\FF_\sigma:\GG_\sigma\to \GG_K$ given by the projections, we have that
   $\kappa^{\FF_\sigma}:F^{\GG_K}\to F^{\GG_\sigma}$ is an equivalence.
  \end{definition}
  We now have all of the relevant definitions to understand the categorical Galois theorem of \cite{galois-theories}, as well as the motivation for the ``Galois'' in ``categorical Galois theorem.''
     \begin{theorem*}[7.5.2 of \cite{galois-theories}]
       \label{thm:1-galois-theories}
       \label{thm:classical-galois-theorem}
       \label{thm:categorical-galois-theorem}
Given
     \begin{align*}
      \C\in \Cat & & F_1,F_2:\C^{op}\to \CCat \\
      \alpha: F_1\tto F_2  & &\sigma: L \to K \text{ (in $\C$)}
     \end{align*}
     where $\alpha$ is componentwise fully faithful, and $\sigma$ is of $F_2$-effective descent,
     we have an equivalence of categories
     \[
     Split_\alpha (\sigma)\simeq F_1^{\GG_\sigma}
     \]

     \end{theorem*}
     Hopefully this expository section makes it clear how classical Galois theory of fields is an instance of the above categorical result.
     For an explanation of how this also subsumes the covering theory of spaces, we refer the interested reader to \cite[Chapter 6]{galois-theories}.

\section{Background}
\label{sec:background}
\label{ch:quasicategorical-galois-theories}
This section is aimed at professional mathematicians in areas roughly adjacent to higher category theory, but does not assume familiarity with the technical aspects of the theory. 
\subsection{Enhanced Mapping Functors}
\label{sec:cart-fibs}
\label{sec:enhanced-mapping}
\textit{(A detailed treatment of the foundations of this subsection can be found in \cite{ghn-fibrations}. This section is a trivial extension of the aforementioned work, but the results added are crucial to the main proof, and thus to the thesis.)}

The \infcat $\iCat$ has both mapping spaces $\Map(\C,\C')$ and mapping $\infty$-categories $\Fun(\C,\C')$ for all pairs of categories $\C,\C'$. Furthermore, treating these as functors, the former factors through the latter as in the following diagram:
\[\begin{tikzcd}
{\iCat\op \times \iCat} && \iCat \\
\\
&& \iSpaces
\arrow["{\Fun}", from=1-1, to=1-3]
\arrow["{(\bullet)^{core}}", from=1-3, to=3-3]
\arrow[""{name=0, anchor=center, inner sep=0}, "{\Map_{\iCat}}"', from=1-1, to=3-3]
\arrow["{(i)}", shorten <=7pt, Rightarrow, 2tail reversed, from=0, to=1-3]
\end{tikzcd}\]
where the arrow labelled (i) is an equivalance in $\Fun(\iCat\op \times \iCat,\iSpaces)$.
Here we are interested in categories other than $\iCat$ which are $\iCat$-enriched in this sense.
\begin{definition}[from \cite{ghn-fibrations}]
	\label{def:enhanced-mapping}
A category $\C$ \de{has an enhanced mapping functor} if the mapping space functor factors up to some equivalence, as in the following diagram:
\[\begin{tikzcd}
{\C\op \times \C} && \iCat \\
\\
&& \iSpaces
\arrow["{\MAP_\C}", from=1-1, to=1-3]
\arrow["{(\bullet)^{core}}", from=1-3, to=3-3]
\arrow[""{name=0, anchor=center, inner sep=0}, "{\Map_\C}"', from=1-1, to=3-3]
\arrow["{(i)}", shorten <=7pt, Rightarrow, 2tail reversed, from=0, to=1-3]
\end{tikzcd}\]
where the arrow labelled (i) is an equivalance in $\Fun(\C\op \times \C,\iSpaces)$.
For $C,C'\in \C$ we let $\vMAP{C,C'}_\C$ denote $\MAP_\C(C,C')$ (or simply $\vMAP{C,C'}$ if the context is clear).
\end{definition}
The following result is central to the framework presented here, as it allows us to work with the collection of functors $\C\op \to \iCat$ as a bicategory.

\begin{lemma}[Lemma 6.4 of \cite{ghn-fibrations}]
\label{lem:functor-category-enhanced-mapping}
If $\D$ has an enhanced mapping functor, then $\D^\C$ has an enhanced mapping functor given by the composition of the following diagram:
\[
\begin{tikzcd}
{\Fun(\C,\D)\op\times \Fun(\C,\D)} && {\Fun(\C\op \times \C,\D\op \times \D)} \\
\\
{\Fun(\Tw(\C)\op, \D\op \times \D)} && {\Fun(\Tw(\C)\op, \iCat)} \\
\\
\iCat
\arrow["\sim", from=1-1, to=1-3]
\arrow["{p_{\Tw(\C)\op}^*}", from=1-3, to=3-1]
\arrow["{(\MAP_\D)_*}", from=3-1, to=3-3]
\arrow["lim", from=3-3, to=5-1]
\end{tikzcd}\]
where $p_{\Tw(\C)\op} : \Tw(\C)\op\to \C\op \times \C$ is the projection from the twisted arrow category.
\end{lemma}
Wherever we have a notion of mapping object, we get a corresponding notion of orthogonality.
In our case, the following notion turns out to have a theory very similar to the notion of orthogonality coming from $\Map(-,-)$.
Thus, it serves to not only simplify the proof of our main result, but to suggest a model-independent approach.
\begin{definition}
	Let $\D$ be a category with an enhanced mapping functor in the sense of Definition~\ref{def:enhanced-mapping}.
  Given arrows $\alpha: F_1 \to F_2$ and $\gamma : H_1\to H_2$ in $\D$, we say $\alpha$ is \de{enhanced-orthogonal} to $\gamma$, denoted $\alpha \porth \gamma$ if the following square is a limit in $\iCat$:
\[\begin{tikzcd}
{\vMAP{F_2, H_1}_\D} & {\vMAP{F_2, H_2}_\D} \\
{\vMAP{F_1, H_1}_\D} & {\vMAP{F_1, H_2}_\D} \\
\arrow[from=1-1, to=1-2]
\arrow[from=1-1, to=2-1]
\arrow[from=1-2, to=2-2]
\arrow[from=2-1, to=2-2]
\end{tikzcd}\]
\end{definition}
While the above makes sense for any category $\D$ with an enhanced mapping functor, we are only interested in a specific case.
Specifically, applying Lemma~\ref{lem:functor-category-enhanced-mapping} to the case $\D=\iCat$ gives us that $\Fun(\C,\iCat)$ has an enhanced mapping functor for any category $\C$.\footnote{Since this is our context of interest, we will let $\vMAP{\, , \,}$ implicitly denote $\vMAP{\, ,\,}_{\Fun(\C, \iCat)}$ for the remainder of this chapter, where $\C$ is understood from context.} Here we seek to understand enhanced orthogonality with respect to $\Fun(\C,\iCat)$.

\begin{lemma}
\label{lem:porth-characterizations}
\label{lem:porth-characterization}
  Given functors $F_1,F_2,H_1,H_2 : \C \to \iCat$ and natural transformations $\alpha: F_1 \tto F_2$ and $\gamma : H_1\tto H_2$, the following are equivalent:
  \begin{enumerate}
  \item[(1)] $\alpha \porth \gamma$ (as arrows in $\Fun(\C, \iCat)$)
    \item[(2)] for every $f:x\to y\in \C$ $\alpha_x \porth \gamma_{y}$ (as arrows in $\iCat$).
  \end{enumerate}
\end{lemma}
\begin{proof}

	We have the following diagram, where the top equivalence is condition (1).
 \newsavebox{\mappingpb}
 \sbox{\mappingpb}{%
   \(
   lim \left( 
	   \begin{tikzcd}[ampersand replacement=\&]
       \& \vMAP{F_2,H_2}\\
      \vMAP{F_1,H_1} \& \vMAP{F_1, H_2}
 \arrow[ from=1-2, to=2-2]
 \arrow[ from=2-1, to=2-2]
   \end{tikzcd}
   \right)
   \)
  }

 \newsavebox{\mappingpblevels}
 \sbox{\mappingpblevels}{%
   \(
   lim \left( 
	   \begin{tikzcd}[ampersand replacement=\&]
       \& \vMAP{F_2(x),H_2(y)}\\
      \vMAP{F_1(x),H_1(y)} \& \vMAP{F_1(x), H_2(y)}
 \arrow[ from=1-2, to=2-2]
 \arrow[ from=2-1, to=2-2]
   \end{tikzcd}
   \right)
   \)
  }
 \newsavebox{\mappingpbexplicit}
 \sbox{\mappingpbexplicit}{%
   \(
   lim \left( \begin{tikzcd}[ampersand replacement=\&]
       \& (\Fun \circ (F_2 \times H_2)\circ p_{\Tw(\C)\op})\\
       (\Fun \circ (F_1 \times H_1)\circ p_{\Tw(\C)\op})\& (\Fun \circ (F_1 \times H_2)\circ p_{\Tw(\C)\op})
 \arrow[ from=1-2, to=2-2]
 \arrow[ from=2-1, to=2-2]
   \end{tikzcd}
   \right)
   \)
  }
 \newsavebox{\mappingpbexplicitlevels}
 \sbox{\mappingpbexplicitlevels}{%
   \(
   lim \left( \begin{tikzcd}[ampersand replacement=\&]
       \& lim_{\Tw(\C)\op}(\Fun \circ (F_2 \times H_2)\circ p_{\Tw(\C)\op})\\
       lim_{\Tw(\C)\op}(\Fun \circ (F_1 \times H_1)\circ p_{\Tw(\C)\op})\& lim_{\Tw(\C)\op}(\Fun \circ (F_1 \times H_2)\circ p_{\Tw(\C)\op})
 \arrow[ from=1-2, to=2-2]
 \arrow[ from=2-1, to=2-2]
   \end{tikzcd}
   \right)
   \)
  }
 \[\begin{tikzcd}
 {\vMAP{F_2,H_1}} & {\usebox{\mappingpb}} \\
 {lim_{\Tw(\C)\op}(\Fun \circ (F_2\times H_1)\circ p_{\Tw(\C)\op})} 
 \arrow["\simeq_{(1)}"{description}, draw=none, from=1-1, to=1-2]
 \arrow["\simeq"{description}, draw=none, from=1-1, to=2-1]
 \end{tikzcd}\]
 Applying the explicit formula for $\vMAP{\bullet,\bullet}$ to the top-right corner and commuting limits then yields
 \begin{align*}
   &\usebox{\mappingpb}\\
   &\simeq  \usebox{\mappingpbexplicitlevels}\\
   &\simeq lim_{\Tw(\C)\op}\left( \usebox{\mappingpbexplicit}\right)
 \end{align*}
Finally, since limits in $\iCat$ are computed pointwise, condition (1) is thus equivalent to the condition that for every $f:x \to y$ in $\C$, we have
   \[\begin{tikzcd}
 {\vMAP{F_2(x), H_1(y)}} & {\usebox{\mappingpblevels}}
 \arrow["\simeq"{description}, draw=none, from=1-1, to=1-2]
 \end{tikzcd}\]
and this is precisely condition (2).
\end{proof}

One can summarize the above by a corollary, which is sufficient for our purposes: ``enhanced orthogonality in $\Fun(\C, \iCat)$ can be checked pointwise''.
This becomes much more useful once we have results about enhanced orthogonality in $\iCat$. 
Thus, we move on to discuss a specific refinement of the (essentially surjective, fully faithful) orthogonal factorization from ordinary category theory.
\subsection{Foundations on (e.s. , f.f.) OFS}
\textit{(A detailed treatment of this section can be found in \cite{kerodon}.)}
\label{sec:ofs}
We begin with the orthogonal factorization system we will focus on in the present work. Namely, the left class will consist of essentially surjective functors, and the right of fully faithful functors, which we now define.
\begin{definition}[\cite{kerodon} 4.6.2.11]
  We say a functor $F: \C \to \D$ of \nifty{1}-categories, is \de{essentially surjective} if and only if the induced set-function $\pi_0(\core{F}):\pi_0(\core{\C})\to \pi_0(\core{\D})$ is surjective.
\end{definition}
\begin{definition}[\cite{kerodon} 4.6.2.1]
  We say a functor $F: \C \to \D$ of \nifty{1}-categories, is \de{fully faithful} if and only if for every pair of objects $X,Y\in \C$, the induced map
  \[
    \iMap_\C (X, Y) \to \iMap_\D (F(X),F(Y))
  \]
  is a homotopy equivalence of Kan complexes.
\end{definition}

Although the present work only utilizes the orthogonality of essentially surjective and fully faithful functors, we include here a larger collection of orthogonal factorization systems which includes these two notions. This is to emphasize the role our specific orthogonal factorization systems plays in the larger project intended to investigate a more general framework.
\begin{definition}[\cite{kerodon} 4.8.5.1]
  We say a functor $F: \C \to \D$ of \nifty{1}-categories, is \de{full} if and only if for every pair of objects $X,Y\in \C$, the induced map
  \[
    \iMap_\C (X, Y) \to \iMap_\D (F(X),F(Y))
  \]
  is surjective on connected components.
\end{definition}
\begin{definition}[\cite{kerodon} 4.8.5.2]
  We say a functor $F: \C \to \D$ of \nifty{1}-categories, is \de{faithful} if and only if for every pair of objects $X,Y\in \C$, the induced map
  \[
    \iMap_\C (X, Y) \to \iMap_\D (F(X),F(Y))
  \]
  is a homotopy equivalence to some summand of $\iMap_\D (F(X),F(Y))$.
\end{definition}
\begin{definition}[\cite{kerodon} 4.8.5.10]

  A functor $F: \C \to \D$ (in $\iCat$) is \de{$n$-full} if:
  \begin{enumerate}
    \item[$n=0$] it is essentially surjective
    \item[$n=1$] it is full
    \item[$n\geq2$] For every $u:X\to Y$ in $\C$, we have 
      $\pi_{n-2}(\iMap_\C (X,Y),u)\to \pi_{n-2}(\iMap_\C (F(X),F(Y)),F(u))$ is injective and
      $\pi_{n-1}(\iMap_\C (X,Y),u)\to \pi_{n-1}(\iMap_\C (F(X),F(Y)),F(u))$ is surjective.
  \end{enumerate}
\end{definition}

\begin{definition}[\cite{kerodon} 4.8.7.1]
  
  For $n\geq 0$, a functor $F: \C \to \D$ (in $\iCat$) is \de{categorically $n$-connective} if it is $m$-full for all $0\leq m \leq n$.
  For $n< 0$, every functor $F: \C \to \D$ (in $\iCat$) is \de{categorically $n$-connective}.
\end{definition}
\begin{definition}[\cite{kerodon} 4.8.7.1]
  For $n\geq -2$, a functor $F: \C \to \D$ (in $\iCat$) is \de{essentially $n$-categorical} if it is $m$-full for all $n+2 \leq m $.
  For $n< -2$, a functor $F: \C \to \D$ (in $\iCat$) is \de{essentially $n$-categorical} if it is a categorical equivalence.
  Let $\essncat{n}$ denote the class of essentially $n$-categorical functors and $\catnconn{n}$ denote the class of categorically $n$-connective functors.
\end{definition}
Given that the proof of Theorem~\ref{thm:categorical-galois-theorem} was essentially a formal consequence of the orthogonality of e.s. functors and f.f. functors, the result below provides a promising generalization of the (e.s., f.f.) orthogonal factorization system.
\begin{corollary}[\cite{kerodon} 4.8.7.18] 
   \label{prop:es-ff-orthogonality}

   For any integer $n$, if $F_1:\A\to \B$ is in $\catnconn{(n+1)}$ and $F_2: \C \to \D$ is in $\essncat{n}$ we have that the square below is a categorical pullback.
\[\begin{tikzcd}
	{\iFun(\B,\C)} & {\iFun(\A,\C)} \\
	{\iFun(\B,\D)} & {\iFun(\A,\D)}
	\arrow["{ \bullet \circ F_2}", from=1-1, to=1-2]
	\arrow["{F_2\circ \bullet}"', from=1-1, to=2-1]
	\arrow["{F_2\circ \bullet}", from=1-2, to=2-2]
	\arrow["{ \bullet \circ F_2}"', from=2-1, to=2-2]
\end{tikzcd}\]
\end{corollary}
The classes above produce a tower of (enhanced) orthogonal factorization systems with (e.s., f.f.) at the bottom. 
\begin{prop}[\cite{kerodon} 4.8.9.1]
  Let $F :\C \to \D$ be a functor and $n$ and integer. Then the following are equivalent:
  \begin{enumerate}
    \item $F$ is in $\catnconn{(n+1)}$
    \item $F \porth \essncat{n}$
  \end{enumerate}
  
\end{prop}
\begin{xample}
Being an essentially surjective functor is equivalent to being orthogonal to all fully faithful functors by the above with $n=-1$.
\end{xample}
Note that the above is enough for the proof of Theorem~\ref{thm:quasicategorical-galois-theorem}. However, future work relies on these being (enhanced) orthogonal factorization systems. This is proven in \cite{kerodon}. To summarize, the above classes arise as consequences of an enriched small object argument, such as the one proven generally in \cite{cat-homotopy-theory} applied to the finite set of maps in \cite[4.8.7.14]{kerodon}. 
\subsection{Segal Groupoid Objects}
\label{sec:groupoids}
This section is essentially a conceptual merger of Section 6.1.2 of \cite{lurie-htt}, and Section 7.3 of \cite{galois-theories}.
In this section and beyond we refer to use `groupoid object' to mean `$\infty$-groupoid object'.

\begin{definition}
A \de{simplicial object} in an $\infty$-category $\C$ is a functor $\so{X}:\DDelta\op \to \C$.
\end{definition}
For simplicity, we visually represent these via the $1$-skeletal diagram corresponding to $\Delta$, and omit the extra data in a functor from $\DDelta$.


\begin{definition}
A simplicial object in $\C$ is a \de{internal groupoid (or Segal groupoid object)} if it satisfies the following internalization of the Kan condition.
  The functor $\so{X}$ takes all colimits of the form on the left to a limit, depicted on the right\footnote{Note this is equivalent to condition $4\prime\prime$ of Proposition 6.1.2.6 of \cite{lurie-htt}.}:
\begin{equation}
\tag{\de{Kan}}
  \begin{tikzcd}[ampersand replacement=\&]
{[0]} \&\& {[m']} \&\& {\so{X}([n])} \& {\so{X}([m'])} \\
\&\&\& \mapsto \\
{[m]} \&\& {[n]} \&\& {\so{X}([m])} \& {\so{X}([0])}
\arrow[from=1-1, to=3-1]
\arrow["S"', from=3-1, to=3-3]
\arrow["{S'}", from=1-3, to=3-3]
\arrow[from=1-1, to=1-3]
\arrow["\Lrcorner"{anchor=center, pos = 0.125, rotate=180, very near start}, phantom, from=3-3, to=1-1, scale=2]
\arrow[from=1-5, to=3-5]
\arrow[from=3-5, to=3-6]
\arrow[from=1-6, to=3-6]
\arrow[from=1-5, to=1-6]
\arrow["\Lrcorner"{anchor=center, pos=0.125}, draw=none, from=1-5, to=3-6]
\end{tikzcd}
\end{equation}

\end{definition}

\begin{definition}
We denote by \de{$\Grpd(\C)$} the full subcategory of $\Fun(\DDelta\op, \C)$ consisting of internal groupoids.
\end{definition}
\begin{definition}
Let $\C$ have finite limits.
Given a map $\lambda: S \to R$ in $\C$, we define the \de{kernel groupoid object} $\GG_\lambda$ of $\lambda$ via the following diagram:
\[\begin{tikzcd}
	\cdots & {S\times_R S\times_R S} & {S\times_R S} & S
	\arrow[shift right=3, from=1-1, to=1-2]
	\arrow[shift left=3, from=1-1, to=1-2]
	\arrow[shift left, from=1-1, to=1-2]
	\arrow[shift right, from=1-1, to=1-2]
	\arrow[shift right=2, from=1-2, to=1-3]
	\arrow[shift left=2, from=1-2, to=1-3]
	\arrow[from=1-2, to=1-3]
	\arrow[shift right, from=1-3, to=1-4]
	\arrow[shift left, from=1-3, to=1-4]
\end{tikzcd}\]
  Given $C\in \C$ we let $const_\C(C)$ denote the Segal groupoid object given by the constant functor whose value is $C$ (i.e. $\GG_{id_C}$.
\end{definition}
\begin{definition}
	Following \cite{galois-theories}, we call the composition $\dYon_\C:\C \into \Fun(\C\op , \S)\into \Fun(\C\op, \iCat)$ the \de{discrete Yoneda embedding}.

\end{definition}
The following definition connects internal groupoids in some fixed $\C$, along with relative presheaves over them, to the enhanced mapping space for $\Fun(\C\op,\iCat)$.
\begin{definition}
  Given a simplicial object $\so{X}$ in $\C$, we define $\grpdEmbedding{\so{X}}:\C\op \to \iCat$ as the colimit, in $\Fun(\C\op,\iCat)$, of $\dYon_\C\circ \so{X}$.
\end{definition}
We present the following lemma, and it's proof as a callback to the thesis. The following lemma is a formal consequence of the properties of $\Ho$ being applied to a proof of a lower analogue. This is one instance of the framework presented here interacting nicely with the diagram in Figure~\ref{fig:ncat-context}. Furthermore, we will see the result itself to be crucial in the proof of Theorem~\ref{thm:quasicategorical-galois-theorem}.
\begin{lemma}
  \label{lem:es-grpd}
  For any Segal groupoid object $\GG$ in $\C$, we have that the canonical map\footnote{ coming from the definition of $\grpdEmbedding{\GG}$ as a colimit over a diagram containing $\GG_0$} ${\grpdEmbedding{ const_\C(\GG_0)}\to \grpdEmbedding{\GG}}$ is pointwise essentially surjective.
\end{lemma}
\begin{proof}
  We note first that for $C\in \C$, we have $\grpdEmbedding{\GG}(C)$ is the colimit of the following diagram:
\[\begin{tikzcd}
	\cdots & {\Map_\C(C,\GG_2)} & {\Map_\C(C,\GG_1)} & {\Map_\C(C,\GG_0)}
	\arrow[shift right=3, from=1-1, to=1-2]
	\arrow[shift left=3, from=1-1, to=1-2]
	\arrow[shift right, from=1-1, to=1-2]
	\arrow[shift left, from=1-1, to=1-2]
	\arrow[shift right=2, from=1-2, to=1-3]
	\arrow[shift left=2, from=1-2, to=1-3]
	\arrow[from=1-2, to=1-3]
	\arrow[shift right, from=1-3, to=1-4]
	\arrow[shift left, from=1-3, to=1-4]
\end{tikzcd}\]
where $\Map_\C(C,\GG_0) \isequiv \grpdEmbedding{const_\C (\GG_0)}(C)$.
By definition, we thus check that $\Ho(\bullet)$  of the above diagram has an essentially surjective map into its colimit (since $\Ho$ preserves colimits):
\[\begin{tikzcd}
\cdots & {\ho{\Map_\C(C,\GG_2)}} & {\ho{\Map_\C(C,\GG_1)}} & {\ho{\Map_\C(C,\GG_0)}} & \Colim (\Ho(\Map_\C(C,\GG_{\bullet})))
	\arrow[shift right=3, from=1-1, to=1-2]
	\arrow[shift left=3, from=1-1, to=1-2]
	\arrow[shift right, from=1-1, to=1-2]
	\arrow[shift left, from=1-1, to=1-2]
	\arrow[shift right=2, from=1-2, to=1-3]
	\arrow[shift left=2, from=1-2, to=1-3]
	\arrow[from=1-2, to=1-3]
	\arrow[shift right, from=1-3, to=1-4]
	\arrow[shift left, from=1-3, to=1-4]
	\arrow[from=1-4, to=1-5]
\end{tikzcd}\]
This follows from Proposition 7.3.6 of \cite{galois-theories}, which says that for simplicial objects $\so{H}$ in $\lFun(\cat{C},\Cat)$ (for any $\cat{C}\in \Cat$, in particular for the terminal category) the canonical map $H_0\to \Colim (\so{H})$ is essentially surjective.
\end{proof}

\section{Quasicategorical Galois Theories}
\label{sec:galois-theorem}
\label{sec:quasicategorical-galois-theories}
We now have all that we need for a quasicategorical analogue to Theorem~\ref{thm:classical-galois-theorem}.
It is worth pointing out that up to this point, all we have done is cite/prove results at the foundations of higher category theory.
That is, the main result of this work is surprisingly formal, given its apparent depth and applicability.
We begin with a refinement of Definition~\ref{def:splitting-object}.
\begin{definition}
Suppose we have two functors $F_1,F_2:\C\op \to \iCat$, a natural transformation ${\alpha: F_1 \to F_2}$, and an arrow $\lambda:x \to y$ in $\C$. Then we define the \de{category of splitting objects} $\Split_\alpha(\lambda)$ via the following pullback in $\iCat$:
\[\begin{tikzcd}
{\Split_\alpha(\lambda)} & {\vMAP{\grpdEmbedding{\GG_R},F_2}} \\
{\vMAP{\grpdEmbedding{\GG_S},F_1}} & {\vMAP{\grpdEmbedding{\GG_S},F_2}}
\arrow[from=1-2, to=2-2]
\arrow[from=2-1, to=2-2]
\arrow[from=1-1, to=1-2]
\arrow[from=1-1, to=2-1]
\arrow["\Lrcorner"{anchor=center, pos=0.125}, draw=none, from=1-1, to=2-2]
\end{tikzcd}\]
\end{definition}
It remains to check that our definition of $\Split_\alpha(\sigma)$ is a sensible one with respect to that of $\lSplit_\alpha(\sigma)$.
\begin{remark}
  
 Suppose we are given a $1$-category $\cat{C}$, two functors $f_1,f_2:\cat{C}\op \to \CCat$, a natural transformation ${\gamma:f_1 \tto f_2}$, and an arrow $\sigma: s\to r$ in $\cat{C}$.
 Recall that we defined $\lSplit_\alpha(\sigma)$ via the following pseudolimit in $\CCat$:
\[\begin{tikzcd}[ampersand replacement=\&]
	{\lSplit_\alpha (\sigma)} \& {f_2(r)} \\
	{f_1(s)} \& {f_2(s)}
	\arrow["{\alpha_s}"', from=2-1, to=2-2]
	\arrow["{f_2(\sigma)}", from=1-2, to=2-2]
	\arrow[dashed, from=1-1, to=2-1]
	\arrow[dashed, from=1-1, to=1-2]
	\arrow["\lrcorner"{anchor=center, pos=0.125}, draw=none, from=1-1, to=2-2]
\end{tikzcd}\]
 We get, via the Duskin nerve, a category $\C$, two functors $F_1,F_2:\C\op \to \CCat$, a natural transformation ${\Gamma:F_1 \tto F_2}$, and an arrow $\Sigma: S\to R$ in $\C$.
 By Yoneda Lemma\footnote{such as that of \cite{riehl-cosmos} or \cite{nima-yoneda}} for $(\infty,2)$-categories, we have that $\vMAP{\grpdEmbedding{\GG_R},F_1}\isequiv F_1(R)$. Thus, applying (the pullback-preserving) $\Ho$ shows that $\Ho (\Split_\Gamma(\Sigma))\isequiv \lSplit_\gamma(\sigma)$.
\end{remark}
As in  the $1$-categorical case, the above is not inherently useful without the $\lambda$ involved satisfying some notion of descent.
\begin{definition}
   Let $F:\C\op \to \iCat$ and consider a map $\lambda:\so{X} \to \so{Y}$ of simplicial objects in $\C$.
  We say $\lambda$ satisfies \de{$F$-effective descent} if the comparison functor $\vMAP{\lambda, F}: \vMAP{\grpdEmbedding{\so{Y}},F}\to \vMAP{\grpdEmbedding{\so{X}},F}$ is a trivial fibration of categories. 
\end{definition}

We present the following result refining Theorem~\ref{thm:categorical-galois-theorem}.

\begin{theorem}
  \label{thm:higher-galois}
  \label{thm:quasicategorical-galois-theorem}
  \label{thm:qcat-galois-theorem}
  \label{thm:qcat-galois}
  Let $\C\in \iCat$ be a category with finite limits. 
  Suppose we have $F_1,F_2 :\C\op \to \iCat$ with a pointwise fully faithful natural transformation $\alpha:F_1\tto F_2$, and we have $\sigma: S \to R$ in $\C$ satisfying $F_2$-effective descent.
  Then $\Split_\alpha(\sigma) \isequiv \vMAP{\grpdEmbedding{\GG_\sigma},F_1}$.
\end{theorem}
\begin{proof}
  By the following diagram, 
\[\begin{tikzcd}
{\Split_\alpha(\sigma)} & {\vMAP{\grpdEmbedding{\GG_R},F_2}} & {\vMAP{\grpdEmbedding{\GG_\sigma},F_2}} \\
{\vMAP{\grpdEmbedding{\GG_S},F_1}} & {\vMAP{\grpdEmbedding{\GG_S},F_2}}
\arrow["\simeq{(1)}", from=1-2, to=1-3]
\arrow[from=2-1, to=2-2]
\arrow[from=1-2, to=2-2]
\arrow[from=1-1, to=2-1]
\arrow[from=1-1, to=1-2]
\arrow["\Lrcorner"{anchor=center, pos=0.125}, draw=none, from=1-1, to=2-2]
\end{tikzcd}\]
where $\simeq_{(1)}$ holds by $F_2$-effective descent, 
it suffices to show that the square below is a pullback in $\iCat$:
\[\begin{tikzcd}
{\vMAP{\grpdEmbedding{\GG_\sigma}, F_1}} & {\vMAP{\grpdEmbedding{\GG_\sigma}, F_2}} \\
{\vMAP{\grpdEmbedding{\GG_S},F_1}} & {\vMAP{\grpdEmbedding{\GG_S}, F_2}}
\arrow[from=1-1, to=1-2]
\arrow[from=1-2, to=2-2]
\arrow[from=1-1, to=2-1]
\arrow[from=2-1, to=2-2]
\end{tikzcd}\]
By definition, this amounts to showing that $ \grpdEmbedding{\gamma} \porth \alpha $, where $\grpdEmbedding{\gamma}:\grpdEmbedding{\GG_S}\to \grpdEmbedding{\GG_\sigma}$ is the map $\grpdEmbedding{ const_\C(S)}\to \grpdEmbedding{\GG_\sigma}$.
By Lemma~\ref{lem:porth-characterization}, this is equivalent to showing that for every $f:x\to y$ in $\C$ we have $ \grpdEmbedding{\gamma}(x)\porth \alpha(y) $.
Proposition~\ref{prop:es-ff-orthogonality} then lets us reduce to showing that $\grpdEmbedding{\gamma}$ is pointwise essentially surjective, which holds by Lemma~\ref{lem:es-grpd}.
\end{proof}


\subsection{Limitations and Future Directions}
\label{sec:conclusion}
We first note that we make no attempt here to connect the notion of descent used here to other variants, such as monadic descent, or the homotopical descent of \cite{homotopical-descent}.

for future directions of the project here established, it is important to note that however important Theorem~\ref{thm:qcat-galois} may be, it is the method of proof which we aim to emphasize here. It was written with an eye toward identifying a more general framework for such results, which we hope is clear (and more accessible) from the exposition here.

The original motivation for generalizing the ordinary categorical Galois theorem of \cite{galois-theories} was to refine the work of \cite{joyal-tierney} on classifying ordinary sheaf toposes via localic $1$-groupoids. 
Recent work (\cite{abfj-sheaves}) has established that one approach to sheaf toposes in the \nifty{1}-categorical sense is via ``lex modalities'', which are orthogonal factorization systems in which the left class is closed under finite limits. 
Note, \cite{abfj-sheaves} focuses on orthogonal factorization systems whose orthogonality is not one enhanced (or enriched) over $\iCat$ but over $\iSpc$. Thus, the following results are promising toward the aim of finding the appropriate level of generality for the proof of Theorem~\ref{thm:quasicategorical-galois-theorem}.
\begin{corollary}[\cite{kerodon} 4.8.5.29]
  For every integer $n$, the class of $n$-full functors is closed under pullbacks along arbitrary functors.
\end{corollary}
\begin{remark}
This immediately implies that $\essncat{n}$  and $\catnconn{(n+1)}$ are closed under pullbacks for all integer $n$.
\end{remark}
With the above result in mind, we expect that:
\begin{enumerate}
\item each of the orthogonal factorization systems, $(\catnconn{(n+1)},\essncat{n})$, will offer something like an ``$(\infty, 2)$-site'' for the ``$(\infty,2)$-topos'' of $(\infty,1)$-categories.\footnote{ (although at the time of writing there seems to be no agreement on the definitions of any of the quoted terms). }
\item we can use this site to induce one on the functor category $\iFun(\C\op, \iCat)$ which seems to give reasonable notion of a sheafy-localization of a $(\infty, 2)$-presheaf category, and
\item the quasicategorical Galois theorem is a consequence of a generalization to an internal $(\infty,2)$-topos analogue. 
\end{enumerate}
Investigating the above is the purpose of the project initiated in the current work.
The trajectory of Section~\ref{sec:ordinary-galois} ought to extrapolate to a more general machine for producing suitable contexts for the desired higher categorical Galois theorem suggested above. Thus, we aim to provide a result which (re)produces instances  of Galois-type results in higher categorical contexts as a formal consequence of one general result being applied to (nice)\footnote{i.e. a suitable generalization of the notion of stably frobenius ajunctions.} adjunctions.
For some specific examples, we offer the following curation of some adjunctions likely to reproduce known Galois-type results as instances of the above machinery.

The first  comes from work of \cite{galois-stable-homotopy} on Galois theory of stable homotopy theories.
\begin{thm*}[Theorem 5.3.6 of \cite{galois-stable-homotopy}]
  The collection of `Galois Categories' (written $\cat{GalCat}$) and that of profinite $1$-groupoids (written $\cat{Pro}(\cat{Gpd}_{fin})$) are contravariantly equivalent (as $2$-categories).
\[\begin{tikzcd}[ampersand replacement=\&]
\cat{GalCat}\op \& \cat{Pro}(\cat{Gpd}_{fin})
\arrow[""{name=0, anchor=center, inner sep=0},  shift left=2, from=1-1, to=1-2]
\arrow[""{name=1, anchor=center, inner sep=0},  shift left=2, from=1-2, to=1-1]
\arrow["\dashv"{anchor=center, rotate=-90}, draw=none, from=0, to=1]
\end{tikzcd}\]
\end{thm*}

The second comes from replacing profinite groupoids with profinite $\infty$-groupoids.
\begin{thm*}[Theorem E.2.4.1 of \cite{lurie-sag}]
Let $\Prof_\infty$ denote the quasicategory of $\pi$-finite $\infty$-groupoids.
Then $\Prof_\infty$ is a reflective subcategory of $\iGroth$.
\[\begin{tikzcd}[ampersand replacement=\&]
\iGroth \& \Prof_\infty
\arrow[""{name=0, anchor=center, inner sep=0},  shift left=2, from=1-1, to=1-2]
\arrow[""{name=1, anchor=center, inner sep=0},  shift left=2, from=1-2, to=1-1]
\arrow["\dashv"{anchor=center, rotate=-90}, draw=none, from=0, to=1]
\end{tikzcd}\]
\end{thm*}
By restricting the adjunction to the essential image of the inclusion, %
this conjecturally produces a classification of certain structured toposes\footnote{ As defined in \cite{carchedi}.} as internal presheaves over ``profinite-homotopy groupoid objects''.
Finally, we hope that the reader appreciates the extent to which the above list is non-exhaustive, as a feature and not a bug.

\printbibliography[title={References}] 
\end{document}